\documentclass[12pt]{article}

\usepackage[margin=1in]{geometry}
\usepackage{amsmath,amsfonts,graphicx,framed}
\usepackage{amssymb}
\usepackage{amsthm}
\usepackage{fancyhdr}
\usepackage{float}
\usepackage{caption}
\usepackage{comment}
\usepackage[utf8]{inputenc}
\usepackage{authblk}
\usepackage[mathscr]{euscript}
\usepackage{mathtools}
\usepackage{xcolor}

\usepackage[hypertexnames=false]{hyperref}

\theoremstyle{plain}% default
\newtheorem{thm}{Theorem}

\newtheorem{prop}{Proposition}
\theoremstyle{definition}
\newtheorem{defn}{Definition}
\newtheorem{rem}{Remark}
\newtheorem{ex}{Example}

\numberwithin{equation}{section}
\numberwithin{thm}{section}
\numberwithin{prop}{section}
\numberwithin{lemma}{section}
\numberwithin{cor}{section}
\numberwithin{defn}{section}
\numberwithin{rem}{section}
\numberwithin{ex}{section}

\newcommand{\pd}{\partial}

\newcommand{\mbb}{\mathbb}

\newcommand{\ep}{\varepsilon}

\newcommand{\re}{\mbb R}
\newcommand{\al}{\alpha}

\newcommand{\Om}{\Omega}
\newcommand{\eqal}[1]{\begin{equation}\begin{aligned}#1\end{aligned}\end{equation}}

\newcommand{\ov}{\overline}

\begin{document}

\title{Removing singularities for fully nonlinear PDEs}

\author{Ravi Shankar}

\date{\today}

\maketitle

\begin{abstract}
    We show removability of half-line singularities for viscosity solutions of fully nonlinear elliptic PDEs which have %are concave and have a Jacobi inequality
    classical density and a Jacobi inequality.  An example of such a PDE is the Monge-Amp\`ere equation, and the original proof follows from Caffarelli 1990.  Other examples are the minimal surface and special Lagrangian equations.  The present paper's quick doubling proof combines Savin's small perturbation theorem with the Jacobi inequality.  
    The method more generally removes singularities satisfying the single side condition.
    %The method more generally removes singularities contained in conical sets with an acute angle.
    %also removes singularities satisfying the interior box condition.  
\end{abstract}

%branch-type singularities, by analogy with branch-cut.  maybe even branch-cut is OK.

\section{Introduction}

%\subsection{Statement of results}

%\textbf{RESULT 1: Monge-Amp\`ere.} 
\subsection{Monge-Amp\`ere}

%We establish a quick proof that the Monge-Amp\`ere equation
We establish a new proof of a known result that the Monge-Amp\`ere equation
\eqal{
\label{MA}
F(D^2u)=\det D^2u=1,\qquad x\in B_1(0)\subset\re^n
}
does not have viscosity solutions singular inside a half-line.  Denote by $\text{sing}(u)\subset B_1(0)$ the relative complement of those points on which $u$ is twice differentiable.  %It is well known that this is a closed set.%By Savin's small perturbation theorem \cite[Theorem 1.3]{S07}, the singular set is closed.

\begin{thm}
\label{thm:MA}
Let $u\in C(\ov{B_1(0)})$ be a convex viscosity solution of \eqref{MA}.  If $\text{sing}(u)$ is contained in a closed line segment with an endpoint in $B_1(0)$, then  $\text{sing}(u)=\varnothing$.
\end{thm}

\begin{rem}[Singularities of Monge-Amp\`ere]

This result originally follows from the regularity theory of Monge-Amp\`ere.  The singular set is well known by Pogorelov's Hessian estimate \cite{P64,P78} to be the set of points where $u$ fails strict convexity, and is closed in $B_1(0)$. %; see alternatively \cite{C90,S07}.  
It follows from Caffarelli \cite{C90} that $u$ only equals a linear function on sets with no extremal points: subspaces.  For example, the Pogorelov \cite{P78} singular solution $u(x)=|(x_1,x_2)|^{1+\al}f(x_3)$ is singular along a full line.
  %It follows from Caffarelli 1990 \cite{C90} that $u$ solving \eqref{MA} does not degenerate, i.e. equal a linear function, on a set with extremal points.  Consequently, a degenerate set must be a subspace, exemplified by Pogorelov's singular solution.  
  The dimension of degeneration is less than $n/2$, as follows from Caffarelli \cite{C93}, where examples of such dimensions were also produced.  
  It follows from Mooney \cite{M15} that the collection of $k$-dimensional degenerate subspaces has Hausdorff dimension less than $n-k$, and that $\text{sing}(u)$ has sharp upper bound of dimension $n-1$.  In dimension two, Alexandrov showed strict convexity \cite{A42}. Singular set propagation to the boundary was obtained for Gaussian curvature and curvature functions by Smith \cite{S12,S24}.  Namely, the singular set is the union of the convex hull of closed subsets of the domain boundary.  
\end{rem}

\begin{rem}[Improvability]
\label{rem:imp}
    Our result can remove half-line singularities, which seems optimal in dimension three.  However, in higher dimensions, the singularities are codimension two planes.  Can one remove half-plane singularities using a similar proof?  This does not seem clear using the elementary methods of this paper.
    %It is unclear to us whether our elementary methods can improve this result.  Examining the methods, we could at least formulate a stronger conjecture which removes singularities contained in half planes of dimension two.  But in dimension three, the Pogorelov solution $u= |(x_1,x_2)|^{1+\al}f(x_3)$ is singular along the $x_3$ axis, which is certainly contained in the closed half-plane $\{(0,a^2,b):a,b\in\re\}$.%  If singular sets contained in half-planes could be eliminated in any dimension, this would rule out Pogorelov in dimension three.
\end{rem}

%\bigskip
%\noindent
%\textbf{RESULT 2: General PDEs: minimal surface, special Lagrangian equation}
\subsection{General PDEs: minimal surface, special Lagrangian equation}

Our next result extends Theorem \ref{thm:MA} to fully nonlinear equations which are classically dense and have a Jacobi inequality, which includes the minimal surface equation
\eqal{
\label{min}
D\cdot\left(\frac{Du}{\sqrt{1+|Du|^2}}\right)=0,
}
and the special Lagrangian equation
\eqal{
\label{slag}
\text{tr}\arctan D^2u=c\in [(n-2)\pi/2,n\pi/2).
}
We encapsulate these equations and Monge-Amp\`ere in the following definition.

\medskip
\begin{defn}
\label{Fdef}
    Let $F:\text{Sym}(n;\re)\times\re^n\to\re$ be Lipschitz and elliptic $F(M+P,p)\ge F(M,p)$ for $P$ positive definite.  We assume the equation is smooth and strictly elliptic at every point of a smooth solution, or $\pd F(M,p)/\pd M>0$ at all points $(M,p)=(D^2u(x),Du(x))$.  The fully nonlinear equation
    \eqal{
    \label{eqF}
    F(D^2u,Du)=0,\qquad x\in B_1(0)\subset\re^n
    }
    has \textbf{classical density} if for every viscosity solution $u\in C(\ov{B_1(0)})$, there exists a sequence $u_k\in C^\infty(\ov{B_1(0)})$ of smooth solutions to \eqref{eqF} converging uniformly to $u$.

    \smallskip
    It is said to \textbf{have a Jacobi inequality} if there exists a Lipschitz positive function $b(M,p)$ satisfying the Jacobi inequality on all such smooth solutions $u$ in the viscosity sense,
    \eqal{
    \label{Jac}
    F^{ij}b_{ij}&\ge 2F^{ij}b_ib_j/b,\\
    F^{ij}&:=\pd F/\pd M_{ij}(D^2u,Du),\\
    b_i&:=\pd b(D^2u,Du)/\pd x^i,
    %\frac{\pd F}{\pd M_{ij}}\frac{\pd ^2}{\pd x^i\pd x^j}b(D^2u,Du)\ge 2b^{-1}(D^2u,Du)\frac{\pd F}{\pd M_{ij}}\frac{\pd }{\pd x^i}b(D^2u,Du)\frac{\pd }{\pd x^j}b(D^2u,Du),
    }
    and the property that $b(D^2u,Du)\le B$ on a domain $\Omega\subset\subset B_1(0)$ implies higher bounds
    \eqal{
    \label{higher}
    |D^ku|\le C(k,B,\Omega',\Omega,\|u\|_{L^\infty(B_1(0))}),\qquad k=0,1,2,3\dots
    }
    on any subdomain $\Omega'\subset\subset \Omega$ and smooth solution $u$.$\qed$
\end{defn}

\begin{thm}
\label{thm:gen}
    Let $u\in C(\ov{B_1(0)})$ be a viscosity solution of equation \eqref{eqF} having classical density and a Jacobi inequality.  If $\text{sing}(u)$ is contained in a line segment with an endpoint in $B_1(0)$, then $\text{sing}(u)=\varnothing$.
\end{thm}

\begin{rem}[New proof for minimal surface]
    The Jacobi inequality for $b=\sqrt{1+|Du|^2}$ is true for minimal surface \eqref{min}, see e.g. \cite[Eq. (16.49)]{GT77}.  The higher bounds \eqref{higher} follow from the $C^{0,1}\to C^{1,\al}$ local estimate of De-Giorgi Nash.  See e.g. \cite[Theorem 16.8]{GT77} for classical solvability of the Dirichlet problem for smooth solutions, then classical density by the comparison principle for viscosity solutions by smoothly approximating the boundary data.  This gives a new proof of half-line singularity removal; the original proof is by interior regularity, which follows from the gradient estimate of Bombieri-De Giorgi-Miranda \cite{BDM69}, see also Trudinger \cite{T72}.
\end{rem}

\begin{rem}[New proof for special Lagrangian equation]
    The Jacobi inequality for $b=(1+\lambda_{max}(D^2u)^2)^{\ep(n)/2}$ and $n\ge 3$ is true for special Lagrangian equation by the proof of \cite[Proposition 2.1]{WY14}.  The higher bounds \eqref{higher} follow follow from the local estimates of Evans-Krylov $C^{1,1}\to C^{2,\al}$, using the level set convexity in \cite{Y06}.  The classical density follows from the classical solvability in Caffarelli-Nirenberg-Spruck \cite{CNS85};  see also the lecture notes of Yuan \cite{Y10} and the ones \cite{Y04,B20} for the comparison principle for viscosity solutions.  This new proof of half-line singularity removal comes after the interior regularity; see \cite{WY09,WY14} for Hessian estimates and \cite{WY10} for gradient estimates.
\end{rem}

\begin{rem}[New proof for Monge-Amp\`ere equation]
\label{rem:MA}
    We show that Theorem \ref{thm:MA} follows from Theorem \ref{thm:gen}.  Let us first note the standard facts that operator $F(D^2u)$ can be defined elliptic everywhere by extending $F=\det D^2u-1$ to $-1$ outside the positive definite matrices, and viscosity solutions are defined by touching $u$ from above and below by convex polynomials.  The comparison principle for viscosity solutions, using again the general proof written down in e.g. \cite{Y04,B20}, works verbatim for the Monge-Amp\`ere equation.

    \smallskip
    Next, the Monge-Amp\`ere equation has a Jacobi inequality \cite{Y23} for $b=\det (I+D^2u)^{1/2n}$.  By Evans-Krylov, the higher estimates \eqref{higher} work.  Finally, the classical solvability \cite{CNS84,K84} and comparison principle for Monge-Amp\`ere implies the classical density.  Thus singularity removal Theorem \ref{thm:MA} works.
\end{rem}

\begin{rem}[Other PDEs with Jacobi inequalities]
    There are many other nonlinear PDEs which have Jacobi or Jacobi type inequalities, at least under some \textit{a priori} conditions on the solutions.  To list a few recent ones, there are the Lagrangian mean curvature equations \cite{B21,L23a,Z23,BW24a,BW24b}, generalized special Lagrangian equation \cite{Z24}, the Hessian quotient equation \cite{L23b,L24}, scalar curvature equations \cite{CJZ24,QZ24,W24}, and 2-Hessian equation with dynamic semi-convexity \cite{ShY23}.  %By the Monge-Amp\`ere equation and Hessian quotient equation \cite{L24}, having a Jacobi inequality is not sufficient for estimates and regularity, but it is a good guess.
\end{rem}

%\noindent
%\textbf{RESULT 3: Extensions}
\subsection{Extensions}

Our last result extends the removable singularities from half-lines to sets satisfying the \textbf{single side condition}, which roughly states that the singularity only intersects a single side of some box.  It excludes singularities which intersect two sides of any closed box.  Half-lines and point singularities are basic examples of single side singularities, while full lines and half-planes fail this condition.

\smallskip
We recall that a \textbf{closed box} in $B_1(0)$ is a rotation about the origin of a set of the form 
$$
\{x\in\re^n:a_i\le x_i\le b_i,1\le i\le n\},\qquad a_i<b_i\quad \text{for all i},
$$
which is contained in $B_1(0)$.  The sides of the box are closed pieces of hyperplanes which comprise its boundary.  We also denote $B^0$ by the interior of a box $B$.%  We state the analogue of Theorems \ref{thm:MA} and \ref{thm:gen} for the interior box condition.

%\begin{defn}
    %Given a set $A\subset B_1(0)$ and a closed box $B=$, we say that the pair $(A,B)$ satisfies the \textbf{interior box condition} if $A\cap \pd B$ is contained in the interior of one of $B$'s sides.  Figure XYZ?
%\end{defn}

\smallskip
The theorem states that the restriction of a singularity to a box is removable if it only intersects a single side of the box.

\begin{thm}[Monge-Amp\`ere]
\label{thm:MA_int}
Let $u\in C(\ov{B_1(0)})$ be a convex viscosity solution of \eqref{MA} on domain $B_1(0)\subset\re^n$.  If there exists a closed box $B$ in $B_1(0)$ such that $\text{sing}(u)$ intersects only one side of $B$, then $\text{sing}(u)\cap B^0=\varnothing$.
\end{thm}

\begin{thm}[General PDEs]
\label{thm:genA}
    Let $u\in C(\ov {B_1(0)})$ be a viscosity solution of \eqref{eqF} having classical density and Jacobi inequality \eqref{Jac}.  If there exists a closed box $B$ in $B_1(0)$ such that $\text{sing}(u)$ intersects only one side of $B$, then $\text{sing}(u)\cap B^0=\varnothing$.
\end{thm}

\begin{ex}[Keyhole singularities]
    By iterating this result, we can remove singularities contained in a keyhole, i.e. a pyramid inside $\Omega$ for which a single vertex intersects the boundary.  %, see Figure XYZ.  
    Moreover, there are removable singularities not contained in a keyhole.  One can remove any singularity whose closure intersects the boundary only once.%, see Figure ABC.
\end{ex}

\subsection{Discussion}
%PDEs such as Monge-Amp\`ere \eqref{MA} have singular viscosity solutions, so the regularity of solutions requires a classification.  
%The question of singularity removal is complementary to partial regularity.  The former seeks to find the largest class of inadmissible singularity, while the latter question is to find the smallest admissible singularity.
%\smallskip
Our results show that the singularities of solutions of a PDE must obey geometric constraints if the PDE has a Jacobi inequality.  Moreover, the Pogorelov singular solution of the Monge-Amp\`ere equation \eqref{MA} shows that the geometric condition in Theorem \ref{thm:MA_int} is rather sharp.  Whether it can be improved is unclear to us; see Remark \ref{rem:imp}.  

\smallskip
An earlier, standard result is the removability of point singularities for $C^\infty$-solvable equations, i.e. those PDEs which have classical solutions to Dirichlet problems with smooth boundary data, such as concave PDEs as in \cite{CNS85}.  A viscosity solution has a point singularity at $p\in\Omega$ if $u\in C^\infty(\ov\Omega/\setminus \{p\})\cap C^0(\ov\Omega)$.  The result states that if a $C^\infty$-solvable equation has a solution with a point singularity, then it is smooth.  The proof is by the comparison principle.%  The boundary data is smooth, 

%An earlier, standard result is the removability of point singularities for $C^\infty$-solvable equations, i.e. those PDEs which have classical solutions to Dirichlet problems with smooth boundary data, as in \cite{CNS85}. Such equations have global \textit{a priori} estimates which control the derivatives by their boundary values.  The boundary data for solutions with point singularities is smooth, so the estimates rule out point singularities.  More generally, the estimates rule out singularities with a positive distance to the boundary.  %  More precisely, the uniqueness of the comparison principle combined with 

\smallskip
In contrast, removing singularities which intersect the boundary is not well understood.  The global estimates break if the boundary values can be singular.  Our results show that if the PDE has a Jacobi inequality, the singularities can be removed even if the boundary values are singular.  The only requirement is that the intersection of the singularity with the boundary be sufficiently small.  For example, a half-line intersects the boundary only at a point.  More general removable singularities satisfy the single side condition, which enforces a weak transversality of the singularity to the boundary.

\smallskip
Other results on the partial regularity for fully nonlinear PDEs include sharp Hausdorff dimension upper bounds of $n-\ep$ in \cite{ASS12} for uniformly elliptic PDEs, and $n-1$ by Mooney \cite{M15} for the Monge-Amp\`ere equation.  Applying Savin \cite{S07} to PDEs with an Alexandrov theorem, such as for convex functions or k-convex functions and solutions, see \cite{CT05} and \cite[Remark 5.2]{ShY23}, shows that the singular set is closed and Lebesgue measure zero.

%Our results are in the same vein as the well known removability of point singularities, or more generally singularities at a positive distance from the boundary, for viscosity solutions of any classically solvable PDE.  That standard fact follows from classically solving the Dirichlet problem on a ball containing the singular point, for which the boundary data is smooth.  Applying the uniqueness granted by the comparison principle shows that the singular point is absent.

%\medskip
%Removing singularities which can intersect the boundary is not well understood, outside of rare situations where full interior regularity is established.  The main obstruction is in finding suitable \textit{a priori} estimates.  The smooth solvabilty in \cite{CNS85} relies on global estimates for the derivatives in terms of their boundary traces.  However, the global estimates fail when we allow singularities which intersect the boundary, such as half-lines.  The interior box condition studied here, which includes half-lines, allows the singularities to intersect the boundary, provided their intersection with the boundary is sufficiently small.

\subsection{Outline of the doubling proof}
Our main observation is a doubling inequality \eqref{doub} which controls $b(D^2u,Du)$ inside $[-1,0]\times[-1,1]$, for example, by its values on three sides of this rectangle.  

\smallskip
A doubling inequality can be thought of as a maximum principle that fails on a small set.  It controls the quantity of interest by its values on that set.  Given a Jacobi inequality, or equivalently a positive superharmonic function, a doubling inequality can often be proved.  One proof is by a comparison argument with a Green's function, as in \cite[Lemma 2.6, pg 4056]{SY23}.  Another proof is by Korevaar type arguments \cite{S24,SY24}.  An integral doubling inequality was shown early by Trudinger \cite[pg 70]{T80} using the ABP estimate.  Using an almost-Jacobi inequality, for which a superharmonic function is not available, the doubling inequality of Qiu \cite{Q24} was extended in \cite{ShY23} using the comparison function of Guan-Qiu \cite{GQ19}.

%\smallskip
%In general, combining a Jacobi inequality with a Green's function cutoff yields a doubling inequality, which controls $b(D^2u,Du)$ by its values near the Green's function pole; see e.g.  for the Green's function comparison argument.  %The doubling inequality \cite{ShY23}, using an approximate Green's function, arises in \cite{T80} from a Korevaar type argument \cite{K87}, and that \cite{ShY23} from an almost-Jacobi inequality.  

\smallskip
Our approach follows \cite{SY23} and compares with a ``1D Green's function" $(1/|x_1|-1)_+$.  In \cite{SY23}, the Green's function vanishes on the boundary, so the defective set is near the pole of the Green's function.  But in our situation, the Green's function is nonvanishing near $x_2=\pm 1$ as well as the pole $x_1=0$.  The doubling inequality only controls $b(D^2u,Du)$ in terms of its values on these sides of a box.  But provided the singular set, such as a half-line, satisfies the single side condition, it follows that $b(D^2u,Du)$ is bounded on the three bad sides, and we obtain a global estimate.  
%But in our situation, the ``1D Green's function" is defective near the transverse boundaries $x_2=\pm 1$ as well as the pole $x_1=0$.  The inequality only controls $b(D^2u,Du)$ in terms of its values on three sides of a box.  Provided the singular set, such as a half-line, satisfies the single side condition, it follows that $b(D^2u,Du)$ is controlled on the three bad sides, and we obtain a global estimate.  

\medskip
The final ingredient in the proof is Savin's small perturbation theorem \cite[Theorem 1.3]{S07}, which lets us link the \textit{a priori} estimate to the partial regularity assumption away from the singular set.  If we use the smooth approximations from classical density, the smooth solutions will be uniformly bounded on the smooth set.  Applying this boundedness to the \textit{a priori} doubling inequality gives global higher derivative bounds on the smooth approximations.  The limit will therefore be interior smooth.

\section{Proof of Theorems}

We prove Theorem \ref{thm:genA}, since the others follow from this one.

\medskip
\noindent
\textbf{Step 1. Single side singularity.  }Let $u$ be a viscosity solution of \eqref{eqF}, and suppose there exists a closed box $B$ in $B_1(0)$ such that $\text{sing}(u)\cap \pd B$ is contained in the interior of one of the sides of $\pd B$.  After a rotation, we write this box as as a set of the form $\{a_i\le x_i\le b_i\}$, and we suppose $\text{sing}(u)\cap \pd B\subset L^0$ for left side $L:=\{x_1=a_1\}\cap B$.  Since $\text{sing}(u)$ is closed, we can find a slightly smaller box $B'$ satisfying the same hypotheses by replacing $b_1$ with $b_1'=b_1-r$ for some $r>0$.

\smallskip
Now let us denote $R:=\pd B'\setminus L^0$ as the ``right" portion of $\pd B'$ which excludes the left side $L$.  The importance of the single side condition is that closed set $\text{sing}(u)$ has a positive distance from $R$.  Indeed, this follows from closed set $\text{sing}(u)\cap \pd B$ being contained in open set $L^0$.  Thus, $R$ is compactly contained in the smooth, open set $\text{sing}(u)^c=B_1(0)\setminus \text{sing}(u)$.

\medskip
\noindent
\textbf{Step 2. Doubling inequality on smooth approximations. }
By classical density, we find smooth solutions $u_k\in C^\infty(\ov{B_1(0)})$ solving \eqref{eqF} which converge to $u$ uniformly as $k\to\infty$.  Letting $b_k=b(D^2u_k,Du_k)$, we compare with a ``1D Green's function" and form the test function on box $B'$
\eqal{
w(x)=\eta(x_1)b_k(x),\qquad \eta(x_1)=\frac{1}{b_1-x_1}-\frac{1}{b_1-a_1}.
}
Suppose $w$ achieves its maximum value $W>0$ inside $B'$ at $x=p$, or $w(x)\le W$ with equality at $x=p$.  Then $b_k(x)\le W/\eta(x_1)$ with equality at $x=p$, i.e. $W/\eta$ touches $b_k$ from above at $x=p$.  Since the Jacobi inequality \eqref{Jac} is true in the viscosity sense, it follows that, at $x=p$,
\eqal{
F^{11}\left(\frac{1}{\eta}\right)''\ge 2F^{11}\eta'(x_1)^2/\eta^3.
}
Since $F^{11}>0$ by ellipticity on smooth solutions, this simplifies to 
\eqal{
0>-\eta''/\eta^2\ge 0.
}
This contradiction implies the maximum occurs on the boundary.  Recalling right portion $R$ in Step 1, we obtain
\eqal{
\label{doub}
\sup_{B'}\eta(x_1)b_k(x)&= \sup_{\pd B'}\eta(x_1)b_k(x)\\
&\le r^{-1}\sup_{R}b_k(x).
}
Thus we obtain control of $b_k$ by its values on a piece of the boundary.

\medskip
\noindent
\textbf{Step 3. Savin stability of partial regularity on $R$.  }By uniform convergence $u_k\to u$ in $B_1(0)$ and Savin \cite[Theorem 1.3]{S07}, it follows that $u_k\to u$ locally smoothly in the smooth set, or in $C^{2,\al}_{loc}(\text{sing}(u)^c)$.  On the other hand, by Step 1, $R$ is compactly contained in $\text{sing}(u)^c$.  It follows that $Du_k(x)$ and $D^2u_k(x)$ are uniformly bounded on $R$ for all $x$ and $k$.  Therefore, there exists $A>0$ such that, for all $k$,
\eqal{
\sup_{B'}\eta(x_1)b_k(x)\le A.
}

\medskip
\noindent
\textbf{Step 4. Interior regularity on $B^0$.}
Fix $x_0\in (B')^0$.  For any $x\in (B')^0$, we find that $b(D^2u_k(x),Du_k(x))=b_k(x)\le C(x_1,A)$ is locally bounded uniformly in $k$, with blowup only as $x_1\to a_1$.  By the higher local bounds \eqref{higher} and eventual boundedness of $\|u_k\|_{L^\infty(B_1)}$, it follows that $|D^mu_k(x)|\le C(m,A,x_1,\|u\|_\infty)$ is bounded nearby $x_0$.  After a subsequence, we assume $u_k\to u$ in $C^{2,\al}$ nearby $x_0$.  We conclude $u$ is $C^{2,\al}$ nearby $x_0$.  Varying $x_0$, this shows $\text{sing}(u)\cap (B')^0=\varnothing$, hence $\text{sing}(u)\cap B^0=\varnothing$.  This completes the proof.

\section{Appendix}

\medskip
\noindent
\textbf{Savin small perturbation theorem.  }We specialize Savin \cite[Theorem 1.3]{S07} to equations $F(M,p,x)=0$ defined on $\text{Sym}(n;\re)\times\re^n\times B_1(0)$, for Definition \ref{Fdef}.

\begin{prop}[\cite{S07}, Theorem 1.3]
\label{prop:Savin}
Let $F(D^2u,Du,x)=0$ satisfy the following hypotheses:

(i) $F$ is locally Lipschitz and elliptic, or $F(M+N,p,x)\ge F(M,p,x)$ if $N>0$, with strict ellipticity in a neighborhood of $M=p=0$, $\pd F/\pd M>0$,

(ii)  $F(0,0,x)=0$ for all $x$,

(iii) $F\in C^2$ in a neighborhood of $M=p=0$.

\noindent
Then there exists $c_1$ small enough depending on $n,F$ such that if viscosity solution $u$ of $F(D^2u,Du,x)=0$ satisfies flatness $\|u\|_{L^\infty(B_1(0))}\le c_1$, then $u\in C^{2,\alpha}(B_{1/2})$ with $\|u\|_{C^{2,\al}(B_{1/2})}\le 1$.  
\end{prop}

\medskip
\noindent
\textbf{Savin stability of partial regularity.  }Suppose $u$ is a viscosity solution of elliptic PDE $G(D^2u,Du)=0$ satisfying the assumptions of Definition \ref{Fdef}, which is smooth on a subdomain $\Omega\subset B_1(0)$, with smooth solutions $u_k\to u$ uniformly on $B_1(0)$.  Then by Savin, $u_k\to u$ in $C^k_{loc}(\Om)$ for any $k$, i.e. $C^k$ convergence in any compact subset of $\Omega$.  Indeed, let $v_k=u_k-u$ solve the elliptic equation on $\Omega$
$$
F(M,p,x)=G(D^2u(x)+M,Du(x)+p)-G(D^2u(x),Du(x))
$$
which is $C^\infty$ near points $M,p=0$ and satisfies Savin's assumptions.  Since $\|v_k\|_\infty\to 0$ as $k\to\infty$, it follows from Savin and a covering argument on small enough balls that $v_k$ has locally uniform-in-$k$ $C^{2,\al}$ and higher estimates as $k\to\infty$.


\begin{thebibliography}{999}


\bibitem[A42]{A42} Aleksandrov, A.D., Smoothness of a convex surface of bounded Gaussian curvature, \textit{Dokl. Akad. Nauk. SSSR} 36 (1942), 195-199.

\bibitem[ASS12]{ASS12} Armstrong, S. N., Silvestre, L. E., \& Smart, C. K. (2012). Partial regularity of solutions of fully nonlinear, uniformly elliptic equations. \textit{Communications on pure and applied mathematics}, 65(8), 1169-1184.

\bibitem[B20]{B20} Bhattacharya, A. (2020). The Dirichlet problem for Lagrangian mean curvature equation. arXiv preprint arXiv:2005.14420.

\bibitem[B21]{B21} Bhattacharya, A. (2021). Hessian estimates for Lagrangian mean curvature equation. \textit{Calculus of Variations and Partial Differential Equations}, 60(6), 224.

\bibitem[BW24a]{BW24a} Bhattacharya, A., \& Wall, J. (2024). Hessian estimates for the Lagrangian mean curvature flow. \textit{Calculus of Variations and Partial Differential Equations}, 63(8), 201.

\bibitem[BW24b]{BW24b} Bhattacharya, A., \& Wall, J. (2024). A Priori Estimates for Singularities of the Lagrangian Mean Curvature Flow with Supercritical Phase. arXiv preprint arXiv:2407.12756.

\bibitem[BDM69]{BDM69} Bombieri, E., De Giorgi, E. \& Miranda, M. (1969) ``Una maggiorazione a priori relative alle ipersuperfici minimali nonparametriche," \textit{Arch. Ration. Mech. Anal.} 32, 255-269.

\bibitem[C90]{C90} Caffarelli, L. A. (1990). A localization property of viscosity solutions to the Monge-Ampere equation and their strict convexity. \textit{Annals of mathematics}, 131(1), 129-134.

\bibitem[C93]{C93} Caffarelli, L. A. (1993). A note on the degeneracy of convex solutions to Monge Amp\`ere equation. \textit{Communications in partial differential equations}, 18(7-8), 1213-1217.

\bibitem[CC95]{CC95} Caffarelli, L. A., \& Cabr\'{e}, X. (1995). Fully nonlinear elliptic equations (Vol. 43). American Mathematical Soc..

\bibitem[CNS84]{CNS84} Caffarelli, L., Nirenberg, L., \& Spruck, J. (1984). The dirichlet problem for nonlinear second‐order elliptic equations I. Monge‐ampégre equation. \textit{Communications on pure and applied mathematics}, 37(3), 369-402.

\bibitem[CNS85]{CNS85} Caffarelli, L., Nirenberg, L., \& Spruck, J. (1985). The Dirichlet problem for nonlinear second order elliptic equations, III: Functions of the eigenvalues of the Hessian. \textit{Acta Mathematica}, 155, 261-301.

\bibitem[CT05]{CT05} Chaudhuri, N., \& Trudinger, N. S. (2005). An Alexsandrov type theorem for k-convex functions. \textit{Bulletin of the Australian Mathematical Society}, 71(2), 305-314.

\bibitem[CJZ24]{CJZ24} Chen, R., Jian, H., \& Zhou, X. (2024). An Integral Approach to Prescribing Scalar Curvature Equations. arXiv preprint arXiv:2408.14850.

\bibitem[GT77]{GT77} Gilbarg, D., Trudinger, N. S., Gilbarg, D., \& Trudinger, N. S. (1977). Elliptic partial differential equations of second order (Vol. 224, No. 2). Berlin: Springer.

\bibitem[GQ19]{GQ19} Guan, P., \& Qiu, G. (2019). Interior $C^2$ regularity of convex solutions to prescribing scalar curvature equations. \textit{Duke Mathematical Journal}, 168(9), 1641-1663.

\bibitem[K87]{K87} Korevaar, N. J. (1987). A priori interior gradient bounds for solutions to elliptic Weingarten equations. \textit{Annales de l'Institut Henri Poincaré C, Analyse non linéaire}, vol. 4, No. 5, pp. 405-421.

\bibitem[K84]{K84} Krylov, N. V. (1984). Boundedly nonhomogeneous elliptic and parabolic equations in a domain. \textit{Mathematics of the USSR-Izvestiya}, 22(1), 67.

\bibitem[L23a]{L23a} Lu, S. (2023). On the Dirichlet problem for Lagrangian phase equation with critical and supercritical phase. \textit{Discrete and Continuous Dynamical Systems}, 43(7), 2561-2575.

\bibitem[L23b]{L23b} Lu, S. (2023). Interior $ C^ 2$ estimate for Hessian quotient equation in dimension three. arXiv preprint arXiv:2311.05835.

\bibitem[L24]{L24} Lu, S. (2024). Interior $ C^ 2$ estimate for Hessian quotient equation in general dimension. arXiv preprint arXiv:2401.12229.

\bibitem[M15]{M15} Mooney, C. (2015). Partial Regularity for Singular Solutions to the Monge‐Ampère Equation. \textit{Communications on Pure and Applied Mathematics}, 68(6), 1066-1084.

\bibitem[P64]{P64} Pogorelov, A. V., \textit{Monge-Amp\`{e}re equations of elliptic type.} Translated from the first Russian edition by Leo F. Boron with the assistance of Albert L. Rabenstein and Richard C. Bollinger. P. Noordhoff, Ltd., Groningen, 1964.

\bibitem[P78]{P78} Pogorelov, A. V., \emph{The Minkowski
multidimensional problem.} Translated from the Russian by Vladimir Oliker.
Introduction by Louis Nirenberg. Scripta Series in Mathematics. V. H. Winston
\& Sons, Washington, D.C.; Halsted Press [John Wiley \& Sons], New
York-Toronto-London, 1978.

\bibitem[Q24]{Q24} Qiu, G. (2024). Interior Hessian Estimates for $\sigma_2$ Equations in Dimension Three. \textit{Frontiers of Mathematics}, 1-22.

\bibitem[QZ24]{QZ24} Qiu, G., \& Zhou, X. (2024). A priori interior estimates for special Lagrangian curvature equations. arXiv preprint arXiv:2407.15159.

\bibitem[S07]{S07} Savin, O. (2007). Small perturbation solutions for elliptic equations. \textit{Communications in Partial Differential Equations}, 32(4), 557-578.

\bibitem[SY23]{SY23} Savin, O., \& Yu, H. (2023). Contact points with integer frequencies in the thin obstacle problem. \textit{Communications on Pure and Applied Mathematics}, 76(12), 4048-4074.

\bibitem[S12]{S12} Smith G., Global singularity theory for the Gauss curvature equation, \textit{Ensaios Matemáticos}, 28, (2015), 1-114, ISSN 2175-0432
%Smith, G. (2012). Global singularity theory for the Gauss curvature equation. arXiv preprint arXiv:1206.5544.
%https://ensaios.sbm.org.br/wp-content/uploads/sites/8/sites/8/2021/11/EM_28-1.pdf

\bibitem[S24]{S24} Smith, G. (2024). On the propagation of singularities of constant curvature, convex hypersurfaces. arXiv preprint arXiv:2410.09950.

\bibitem[Sh24]{Sh24} Shankar, R. (2024). Hessian estimates for special Lagrangian equation by doubling. arXiv preprint arXiv:2401.01034.

\bibitem[ShY23]{ShY23} Shankar, R., \& Yuan, Y. (2023). Hessian estimates for the sigma-2 equation in dimension four. arXiv preprint arXiv:2305.12587.

\bibitem[SY24]{SY24} Shankar, R., \& Yuan, Y. (2024). Regularity for the Monge–Ampère equation by doubling. \textit{Mathematische Zeitschrift}, 307(2), 34.

\bibitem[T72]{T72} Trudinger, N. S. (1972). A new proof of the interior gradient bound for the minimal surface equation in n dimensions. \textit{Proceedings of the National Academy of Sciences}, 69(4), 821-823.

\bibitem[T80]{T80} Trudinger, N. S. (1980). Local estimates for subsolutions and supersolutions of general second order elliptic quasilinear equations. \textit{Inventiones mathematicae}, 61(1), 67-79.

\bibitem[W24]{W24} Wang, B. (2024). Hypersurfaces of constant scalar curvature in hyperbolic space with prescribed asymptotic boundary at infinity. arXiv preprint arXiv:2408.07656.

\bibitem[WY09]{WY09} Warren, M., \& Yuan, Y. (2009). Hessian estimates for the sigma‐2 equation in dimension 3. \textit{Communications on Pure and Applied Mathematics}, 62(3), 305-321.

\bibitem[WY10]{WY10} Warren, M., \& Yuan, Y. (2010). Hessian and gradient estimates for three dimensional special Lagrangian equations with large phase. \textit{American Journal of Mathematics}, 132(3), 751-770.

\bibitem[WY14]{WY14}  Wang, D., \& Yuan, Y. (2014). Hessian estimates for special Lagrangian equations with critical and supercritical phases in general dimensions. \textit{American Journal of Mathematics}, 136(2), 481-499.

\bibitem[Y04]{Y04} Yu Yuan, Lecture notes on linear and nonlinear elliptic equations.

\bibitem[Y06]{Y06} Yuan, Y. (2006). Global solutions to special Lagrangian equations. \textit{Proceedings of the American Mathematical Society}, 134(5), 1355-1358.

\bibitem[Y10]{Y10} Lecture notes on Dirichlet problem for special Lagrangian equations, a model case.

\bibitem[Y23]{Y23}Yuan, Y. (2023). A monotonicity approach to Pogorelov’s Hessian estimates for Monge-Ampere equation. \textit{Mathematics in Engineering}, 5(2), 1-6.

\bibitem[Z23]{Z23} Zhou, X. (2023). Hessian estimates for Lagrangian mean curvature equation with sharp Lipschitz phase. arXiv preprint arXiv:2311.13867.

\bibitem[Z24]{Z24} Zhou, X. (2024). Notes on generalized special Lagrangian equation. \textit{Calculus of Variations and Partial Differential Equations}, 63(8), 197.

\end{thebibliography}
\end{document}